\documentclass[11pt,english]{article}
\usepackage[utf8]{inputenc}
\usepackage{mystyle}

\usepackage{cleveref}
\def\abhiOpt#1{{}{\textcolor{red}{#1} }}
\def\taoOpt#1{{}{\sf\textcolor{blue}{#1}}}
\def\lianaOpt#1{{}{\sf\textcolor{cerisepink}{#1}}}

\newcommand{\lam}{\lambda}

\DeclareMathOperator{\E}{E}
\DeclareMathOperator{\ex}{ex}

\title{Rainbow Tur\'an number of clique subdivisions}

\author{
	Tao Jiang\thanks{
	Department of Mathematics, 
	Miami University, 
    Oxford, OH 45056, US. 
	Email: \texttt{jiangt}@\texttt{miamioh.edu}. 
	Research supported by National Science Foundation grant DMS-1855542.
	}
	\and
	Abhishek Methuku\thanks{
School of Mathematics,
University of Birmingham,
Birmingham, B15 2TT, UK. 
	Email: \texttt{abhishekmethuku}@\texttt{gmail.com}.
	Research supported by the EPSRC, grant no. EP/S00100X/1.
	}
	\and
	Liana Yepremyan\thanks{ 
	Department of Mathematics, Emory University, Atlanta, Georgia, 30322, U.S.
	Email: \texttt{lyeprem}@\texttt{emory.edu}. Research supported by Marie Sklodowska Curie Global Fellowship, H2020-MSCA-IF-2018:846304.
 }
}

\begin{document}
\maketitle

\begin{abstract}
We show that for any integer $t\geq 2$, every properly edge-coloured graph on $n$ vertices with more than $n^{1+o(1)}$ edges contains a rainbow subdivision of $K_t$. Note that this bound on the number of edges is sharp up to the $o(1)$ error term.  This is a rainbow analogue of some classical results on clique subdivisions and extends some results on rainbow Tur\'an numbers. Our method relies on the framework introduced by Sudakov and Tomon~\cite{sudakov2020extremal} which we adapt to find robust expanders in the coloured setting.
\end{abstract}

\section{Introduction}

For a graph $H$, the Tur\'an number $\ex(n, H)$ is the maximum number of edges in an $n$-vertex graph without a copy of $H$. Keevash, Mubayi, Sudakov and Verstra\"ete \cite{keevash2007rainbow} introduced a rainbow variant of the Tur\'an problem. In an edge-coloured graph, we say that a subgraph is \emph{rainbow} if no two of its edges have the same colour. The rainbow Tur\'an number $\ex^*(n, H)$ is the maximum number of edges in a properly edge-coloured graph on $n$ vertices which does not contain a rainbow copy of $H$. Clearly, $\ex^*(n,H) \ge \ex(n,H)$ for every $n$ and $H$. Among other things, Keevash, Mubayi, Sudakov and Verstra\"ete~\cite{keevash2007rainbow} proved that if $H$ is non-bipartite then $\ex^*(n,H) = (1+o(1)) \ex(n,H)$ where the asymptotic value of $\ex(n,H)$ is known by celebrated Erd\H{o}s–Stone–Simonovits theorem. However, if $H$ is bipartite much less is known.  Keevash, Mubayi, Sudakov and Verstra\"ete~\cite{keevash2007rainbow}  showed that $\ex^*(n, C_{2k}) = \Omega(n^{1 + 1/k})$ for all $k \ge 2$, and conjectured that $\ex^*(n, C_{2k}) = \Theta(n^{1 + 1/k})$. They verified this conjecture for $k \in \{2, 3\}$ and for general $k$, Das, Lee and Sudakov showed that for every $k \ge 2$, $\ex^*(n, C_{2k}) =  O(n^{\frac{(1 + \varepsilon_k) \log k}{k}})$ where $\varepsilon_k \to \infty$ as $k \to \infty$. Very recently, Janzer~\cite{janzer2020rainbow} proved the conjecture by showing that $\ex^*(n, C_{2k}) = O(n^{1+1/k})$.

It is well-known that a graph on $n$ vertices without any cycle has at most $n-1$ edges. It is then natural to ask how many edges a properly edge-coloured $n$-vertex graph can have if it does not contain any rainbow cycle.  Keevash, Mubayi, Sudakov and Verstra\"ete~\cite{keevash2007rainbow} showed that there are graphs with $\Omega(n \log n)$ edges that can be properly coloured with no rainbow cycle by colouring an $m$-dimensional cube $Q$ as follows.  The vertices of $Q$ are subsets of $\{1, 2, \dots, m\}$ and for any $S \subseteq \{1, 2, \dots, m\}$ and any $i \in S$, there is an edge between $S$ and $S \setminus \{i\}$ of colour $i$. To see that $Q$ has no rainbow cycle, consider any cycle $C$ in $Q$ and an edge $e$ from $A$ to $A \setminus \{i\}$ of colour $i$ on $C$. Observe that there is a subpath of $C$ connecting $A \setminus \{i\}$ back to $A$ which does not contain $e$ but contains at least one other edge of colour $i$, as $A$ contains the element $i$. It follows that $C$ is not a rainbow cycle. Moreover, since $2^m = n$, $Q$ has $\frac{1}{2}mn = \frac{1}{2}n \log n$ edges, as desired. 
On the other hand, Keevash, Mubayi, Sudakov and Verstra\"ete  showed an upper bound of $O(n^{4/3})$. Das, Lee and Sudakov~\cite{das2013rainbow} improved this bound as follows. 

\begin{theorem}[Das, Lee, Sudakov~\cite{das2013rainbow}]
\label{DasLS}
If $\eta > 0$ and $n$ is sufficiently large, then any properly edge-coloured $n$-vertex graph with at least $n e^{(\log n)^{\frac{1}{2}+\eta}}$ edges contains a rainbow cycle. 
\end{theorem}

Recently Janzer~\cite{janzer2020rainbow} improved the bound in Theorem~\ref{DasLS}  by showing that a properly edge-coloured graph on $n$ vertices with average degree at least $c(\log n)^4$ contains a rainbow cycle, for some constant $c>0$. The main result of our paper extends Theorem~\ref{DasLS} by establishing an analogous theorem for  the family of rainbow subdivisions of the clique as stated below.

\begin{theorem}\label{thm:main}
For every integer $t \geq 2$ there exists a constant $c>0$ such that for every integer $n\geq t$ if $G$ is a properly edge-coloured graph  on $n$ vertices with at least $ne^{c\sqrt{\log n}}$ edges then $G$ contains a rainbow subdivision of $K_t$, where each edge of $K_t$ is replaced with a path of length at most $1300 \log^2 n$.
\end{theorem}

    Note that the problem we address here is a ``rainbow'' instance of a well studied class of problems that explore edge density conditions forcing the existence of a minor/subdivision of small size.  A well-known result of Mader~\cite{mader1967homomorphieeigenschaften} states that for every integer $t$, there is a constant $c(t)$ such that every graph with $c(t) n$ edges contains a $K_t$-minor. Estimates on $c(t)$ have been improved in \cite{kostochka1984lower, thomason1984extremal, thomason2001extremal}. Fiorini, Joret, Theis and Wood~\cite{fiorini2012small} asked the following natural question: how many edges suffice to guarantee that a graph contains not only a $K_t$-minor, but one which has few vertices? Addressing this question, Shapira and Sudakov~\cite{shapira2015small} showed that for every $\varepsilon > 0$ and $t \ge 3$, there exists $n_0$ such that every graph on $n \ge n_0$ vertices and $(c(t) + \varepsilon)n$ edges contains a $K_t$-minor of order $O_{\varepsilon, t}(\log n\log\log n)$. Montgomery~\cite{montgomery2015} subsequently improved this bound on the order
of the $K_t$-minor to the optimal bound $O_{\varepsilon, t} (\log n)$.
In confirming a conjecture of Erd\H{o}s, it was shown by Kostochka and Pyber~\cite{kostochka1988small} that any graph with $4^{t^2} n^{1+\varepsilon}$
edges contains a $K_t$-subdivision of order $O(t^2 \log t/\varepsilon)$. The bound  on the order of the $K_t$-subdivision was improved to the optimal
$O(t^2/\varepsilon)$ by Jiang \cite{jiang2011minor}. Theorem~\ref{thm:main} is a a rainbow version of these  results. Our proof establishes an interesting connection between rainbow Tur\'an problems and expanders, building on the method of Sudakov and Tomon used in~\cite{sudakov2020extremal}.

\subsection*{An overview of the proof} Our method builds on the method used by Sudakov and Tomon in \cite{sudakov2020extremal} together with some new ideas. We incorporate the minimality notion commonly used in the study of graph minors and adapt the notion of ``expander'' conveniently to our setting (see Definition~\ref{def:expander-new}). We show that most of the edges of a sufficiently dense graph  can be covered by edge-disjoint expanders.  In a properly edge-coloured expander, from any given vertex $v$, we can reach almost all of the other vertices  by a rainbow path of poly-logarithmic length avoiding a given set of vertices and colours. Additionally, the notion of minimality ensures that the set of these ``reachable'' vertices induces most of the edges in the expander (see Lemma~\ref{lem:reachable-new}). This additional feature of the expander is used to show the existence of a common large intersection of reachable vertices for a pair of vertices $x$ and $y$ in a general graph (not necessarily an expander). Eventually we are either able to join any pair of vertices by a rainbow path of poly-logarithmic length avoiding a bounded set of colours and vertices (thus allowing us to build a copy of the desired rainbow $K_t$-subdivision) or find a much denser subgraph (see Lemma \ref{rainbow-connection}). We then complete the proof via a density increment argument as in \cite{sudakov2020extremal}.

\subsection*{Notation and Standard Tools}  Let $G$ be a graph, we denote its vertex set by $V(G)$ and the edge set by $E(G)$, and let $v(G)=|V(G)|$, $e(G)=|E(G)|$. For any subset $S$ of $V(G)$ let $e(S)=e(G[S])$. Similarly, for two disjoint subsets $S$ and $T$ of $V(G)$ write $e(S, T) = e(G[S, T])$. Let $d(G)$ be the average degree of $G$.
    For a properly coloured graph $F$, let $C(F)$ denote the set of colours used on the edges of $F$.
\begin{lemma}[Chernoff bounds, \cite{molloy2013graph}]Given a binomially distributed variable $X\in Bin(n, p)$  for all  $0<a\leq 3/2$ we have
 $$\Prob{|X-\E[X]|\geq a \E[X]}\leq 2e^{-\frac{a^2}{3}\E[X]}$$.
\end{lemma}
\section{Expanders}
 
\begin{definition} \label{def:d-minimal}
We say that $G$  is $d$-\emph{minimal} if $d(G) \geq d$, but $d(F) < d$ for every proper subgraph $F\subseteq G$.
\end{definition}
Note that for any graph $G$ with $d(G)\geq d$, there exists a smallest subgraph $H$ that satisfies $d(H)\geq d$.
Such a subgraph $H$ is a $d$-minimal subgraph of $G$. Note also that, by definition, if $G$ is $d$-minimal, then it is also
$d(G)$-minimal.
We will utilize the minimality notion in two ways.
One is based on the definition, which, roughly speaking, says that any subset of vertices spans at most the expected number of edges. The other way we utilize minimality is that whenever we delete a set of vertices we must lose at least the expected number of edges, as stated in the following lemma. 

\begin{lemma} \label{lem:d-expand}
If $G$ is $d$-minimal, then for every set $S\subseteq V(G)$, the number of edges having at least one endpoint in $S$ is at least $d|S|/2$, that is
$$e(S)+e(S,S^c)\geq \frac{d |S|}{2}.$$
In particular, $\delta(G) \ge \frac{d}{2}$.
\end{lemma}

\begin{proof} Let $S\subseteq V(G)$. Suppose for a contradiction that $e(S)+e(S,S^c) < \frac{d |S|}{2}$. Then 
$$e(S^c) > \frac{d|V|}{2} - \frac{d|S|}{2} \geq \frac{d(|V| - |S|)}{2}, $$
a contradiction to the $d$-minimality of $G$.
\end{proof}
    
\begin{definition}\label{def:expander-new}
        Given $\eps,\lambda \in (0, 1)$ and $d \ge 1$, an $n$-vertex graph $G$ is called a \emph{$(d,\lambda, \eps)$-expander} if $G$ is $d$-minimal, and for every subset $S \subseteq V(G)$ of size at most $(1 - \eps)n$, we have $d(S) \le (1 -\lambda)d$. 
\end{definition}

Note that by the remarks after Definition~\ref{def:d-minimal}, if $G$
is  a $(d,\lambda,\varepsilon)$-expander then it is also a
$(d(G), \lambda, \varepsilon)$-expander. Hence, we will often say
let $G$ be a $(d,\lambda,\varepsilon)$-expander with $d=d(G)$.
    
    \begin{lemma} \label{obs:edge-expansion-new}
        Let $n, d \ge 1$ and let $\lambda, \eps \in (0, 1)$.
        Let $G$ be a $(d,\lambda, \eps)$-expander on $n$ vertices. Then every $S \subseteq V(G)$ with $|S| \le (1 - \eps)n$ satisfies $e(S, S^c) \ge \frac{\lambda d}{2}|S|$.
    \end{lemma}
    
    \begin{proof}
        Let $S \subseteq V(G)$ and let $|S| \leq (1-\eps)n$. Since $G$ is a $(d,\lambda, \eps)$-expander, it is $d$-minimal by definition. So by Lemma~\ref{lem:d-expand}, $e(S) + e(S, S^c) \ge \frac{d|S|}{2}$ and $e(S) = \frac{d(S)}{2}|S| \le \frac{(1 -\lambda)d}{2}|S|$. It follows that $e(S, S^c) \ge \frac{\lambda d}{2}|S|$, as claimed.
    \end{proof}

    \begin{lemma} \label{lem:existence-expanders-new}
        Let $n$ be a positive integer. Let $d,\varepsilon,\lambda$ be positive reals where $d\geq 1, \varepsilon\in (0,1)$
        and $\lambda \leq \frac{\eps}{2\log n}$.
        Let $G$ be a graph on $n$ vertices with average degree $d$. Then $G$ contains a $(d',\lambda, \eps)$-expander, with $d' \ge \frac{d}{2}$. 
    \end{lemma}
    
    \begin{proof}
        Let $G_0=G$. We run the following process for $i\geq 0$:
        \begin{enumerate}[label = \rm(\alph*)]
            \item 
                Let $d_i=d(G_i)$. Let $H_i$ be a $d_i$-minimal subgraph of $G_i$. 
            \item 
                If $H_i$ is a $(d_i,\lambda, \eps)$-expander, stop. Otherwise, there exists $S_{i+1} \subseteq V(H_i)$, with $|S_{i+1}| \le (1 - \eps)|V(H_i)|$, such that $d(S_{i+1}) \ge (1 -\lambda)d_i$. Take $G_{i+1} = G[S_{i+1}]$. 
        \end{enumerate}
        Suppose that the process stopped right after $H_{\ell}$ was defined; so $H_{\ell}$ is a $(d_{\ell},\lambda, \eps)$-expander. Since $|V(G_{i+1})| \le (1 - \eps)|V(G_i)|$ for every $i \ge 1$ and $|V(G_{\ell})| \ge 1$, it follows that $(1 - \eps)^{\ell}n \ge 1$, implying that $\ell \le \frac{-\log n}{\log(1 - \eps)} \le \frac{\log n}{\eps}$. Since $d_{i+1} \ge (1 -\lambda)d_i$ for every $i \ge 1$ and $d_0 \ge d$, we have $d_{\ell} \ge (1 -\lambda)^{\ell} d \ge (1 -\lambda \ell)d \ge (1 - \frac{\lam\log n}{\eps})d \ge \frac{d}{2}$.
    \end{proof}

The next corollary says that we can iterate the previous lemma to cover almost all the edges of any graph by expanders.

\begin{corollary} \label{cor:covering}
Let $n$ be a positive integer. Let $d,\varepsilon,\lambda$ be positive reals where $d\geq 1, \varepsilon\in (0,1)$ and $\lambda \leq \frac{\eps}{2\log n}$.
Let $G$ be a graph on $n$ vertices with average degree $d$.
Then $G$ contains edge-disjoint subgraphs $G_1, G_2, \dots, G_k$ such that  for each $i\in [k]$, $G_i$ is a $(d(G_i),\lambda, \varepsilon)$-expander, where $d(G_i)\geq \frac{\varepsilon d}{2}$, and $e(\cup_{i=1}^k{G_i})\geq (1-\varepsilon) e(G)$.
\end{corollary}
\begin{proof}
Let $\mathcal{F}=\{G_1, G_2, \ldots, G_k\}$ be a maximal collection of edge-disjoint subgraphs such that for each $i\in [k]$, $G_i$ is a $(d(G_i),\lambda,\varepsilon)$-expander, where $d(G_i)\geq \frac{\varepsilon d}{2}$.
Let $G'$ be the subgraph of $G$ consisting of edges not covered by $G_1, G_2, \ldots, G_k$.
If $e(G') \ge \varepsilon e(G)$ then $d(G') \ge \varepsilon d$.
By Lemma~\ref{lem:existence-expanders-new},
 $G'$ contains a $(d_{k+1}, \lambda,\varepsilon)$-expander $G_{k+1}$ with $d(G_{k+1})=d_{k+1}\geq \frac{\varepsilon d}{2}$, contradicting the maximality of $\mathcal{F}$. 
 It follows that $e(G') \le \eps e(G)$, as required.
\end{proof}

Given an edge-coloured graph and some path $P$ in it, we say P \textit{avoids} a colour $c$ if it does not contain any edge of colour $c$. More generally, we  say $P$ avoids a set of colours $C$  if $P$ avoids each colour $c\in C$. Analogously, we say $P$ \textit{avoids} a vertex $v$ if it does not contain  $v$, and it avoids a set of vertices $F$ if it avoids each vertex $v\in F$.
 
    \begin{lemma} \label{lem:reachable-new}
        Let $n, \ell, M \ge 1$ be integers and let $d,\eps, \lambda$
        be reals where $d\geq 1$ and $\varepsilon,\lambda\in (0, 1)$. Suppose that $\ell = \left \lceil\frac{4 \log n}{\lam}\right \rceil$ and $d \ge \frac{4(\ell +2M)}{\lam}$. Suppose that $G$ is a properly edge-coloured $(\lam, d, \eps)$-expander on $n$ vertices, let $x \in V(G)$ and let $F$ be a set of at most $M$ vertices and $C$ be a set of at most $M$ colours. 
        Then  there is a set of vertices $U$ such that each vertex $u\in U$ can be reached from $x$ by a rainbow path of length at most $\ell+1$ avoiding $C$ and $F$, and furthermore,
        
          \begin{enumerate}
      \item [(i)] $|U|\geq (1-\varepsilon)n$, and
      \item [(ii)] $e(G[U])\geq (1-\lambda/2-\varepsilon) e(G)$. 
  \end{enumerate}
    \end{lemma}
    
    \begin{proof}
    By the remarks after Definition~\ref{def:expander-new}, we may assume that $d(G)=d$.
        For each $i \ge 0$, let $U_i$ be the set of vertices reachable from $x$ by a rainbow path of length at most $i$ avoiding $C$ and $F$. For each $u \in U_i$, fix a rainbow path $P(u)$ from $x$ to $u$ of length at most $i$ that avoids $C$ and $F$ .

         \begin{claim}
        \label{edgesexpand}
            For any $i \le \ell$, $e(U_{i+1}) \geq \left(1 - \frac{\lam}{2}\right)\frac{d|U_i|}{2}.$
        \end{claim}
        
        \begin{proof}
        Write $E_i = E(U_i, U_i^c)$, and let $E_i^*$ be the set of edges $uw$, with $u \in U_i$ and $w \in U_i^c$, such that  either the colour of $uw$ appears on $P(u)$ or in $C$ or $w\in F$. 
            Then 
            \begin{equation} \label{eq:E*-bound}
            |E_i^*| \le (i + 2M)|U_i| \le (\ell+2M) |U_i| \le \frac{\lam d}{4}|U_i|.
            \end{equation}
            
            Consider an edge $uw \in E_i \setminus E_i^*$, with $u \in U_i$ and $w \in U_i^c$. Then, by defintion, $P(u)\cup uw$ is a rainbow path from $x$ to $w$ of length at most $i + 1$, avoiding $F$ and $C$, so $w \in U_{i+1}$, implying that $uw \in E(U_i, U_{i+1} \setminus U_i)$. It follows that $E_i \setminus E_i^* \subseteq E(U_i, U_{i+1} \setminus U_i)$. Hence, for any $i \le \ell$,
            \begin{equation*}
            \label{eq:lowerboundeUi+1}
                e(U_{i+1}) 
                \ge e(U_i) + |E_i \setminus E_i^*| 
                = e(U_i) + e(U_i, U_i^c) - |E_i^*| 
                \ge \frac{d|U_i|}{2} - \frac{\lam d|U_i|}{4} 
                = \left(1 - \frac{\lam}{2}\right)\frac{d|U_i|}{2},
            \end{equation*}
            where we used Lemma~\ref{lem:d-expand} and \eqref{eq:E*-bound}. This proves the claim.
        \end{proof}
        
        \begin{claim}
        \label{claim1expand}
            Let $i \le \ell$. If $|U_{i+1}| \le (1 - \eps)n$, then $|U_{i+1}| \ge (1 + \frac{\lam}{2})|U_i|$.
        \end{claim}
        
        \begin{proof}
           Suppose that $|U_{i+1}| \le (1 - \eps)n$. Since $G$ is a $(\lambda,d,\varepsilon)$-expander and $|U_{i+1}| \le (1 - \eps)n$, by definition, we have $e(U_{i+1}) \le \frac{(1 -\lambda)d}{2}|U_{i+1}|$.
            Combining this with the lower bound on $e(U_{i+1})$ obtained from Claim~\ref{edgesexpand}, we get $|U_{i+1}| \ge \frac{1 -\lambda/2}{1 -\lambda}|U_i| \ge (1 + \frac{\lam}{2})|U_i|$.
        \end{proof}
        Suppose that $|U_{\ell}| \le (1 - \eps)n$. By iterating Claim~\ref{claim1expand} we get $|U_{\ell}| \ge (1 + \frac{\lam}{2})^{\ell} \ge e^{\frac{\lam \ell}{4}}\ge e^{\log n} > n$, a contradiction. This implies that $|U_{\ell}| \ge (1 - \eps)n$. Now, by definition, $U_{\ell+1}\supseteq U_{\ell}$, hence $|U_{\ell+1}|\geq (1-\varepsilon)n$ and furthermore, by Claim~\ref{edgesexpand},
        $$    e(U_{\ell+1}) \geq \left(1 - \frac{\lam}{2}\right)\frac{d|U_{\ell}|}{2} \geq (1-\lambda/2)(1-\varepsilon) e(G)\geq (1-\lambda/2-\varepsilon) e(G),$$
        where in the second inequality we used that $|U_{\ell}| \geq (1-\eps)n$ and  $e(G)=dn/2$. 
        So it follows that  $e(U_{\ell+1})\geq (1-\lambda/2 -\varepsilon) e(G)$.
        The lemma holds with $U=U_{\ell+1}$.
    \end{proof}

\section{Proof of Theorem~\ref{thm:main}}

\begin{lemma} \label{lem:edge-to-vertex}
Let $G$ be a graph with $\delta(G)\geq \frac{1}{3}d(G)$ and  $G'$  be a subgraph of $G$. For any $\gamma >0$, if $e(G')\geq (1- \gamma) e(G)$ then $v(G')\geq (1 - 3 \gamma) v(G)$.
\end{lemma}

\begin{proof}
Let $V=V(G)$. Let $S=V\setminus V(G')$. Suppose for contradiction that $|S| > 3 \gamma \cdot v(G)$.
The number of edges in $G$ that have at least one endpoint in $S$ is at least
$$\frac{\delta(G)|S|}{2} \geq \frac{d(G)}{3} \cdot \frac{3 \gamma \cdot v(G)}{2} > \gamma  \cdot e(G).$$ 
This is impossible
since the set of edges of $G$ that are incident with $S$ is disjoint from $E(G')$ and
$e(G')\geq (1- \gamma) e(G)$. 
\end{proof}

The following is our main lemma of the paper.

\begin{lemma} \label{rainbow-connection} 
Let $0<\varepsilon \leq\frac{1}{40}$ be a fixed real. For all sufficiently large integers $n$ the following holds. Let $\lambda = \frac{\varepsilon}{2\log n}$. Let $\ell=\left \lceil\frac{4\log n}{\lambda}\right \rceil$.
Let $K,L$ be positive integers. Let $M \ge L+ \frac{12 K}{\eps}(\ell+2)$.
Let $d$ be a real that satisfies  $d\geq \frac{24(\ell+2M)}{\varepsilon\lambda}$.
Let $G$ be a properly edge-coloured  $(d, \lambda, \varepsilon)$-expander,
with $d=d(G)$. Suppose $G$ contains no subgraph on at most $\frac{n}{K}$ vertices with average degree at least $\frac{\varepsilon d}{6}$. Then for any two vertices $x,y\in V(G)$,
any set $C$ of at most $L$ colours and any set $F$ of at most $L$ vertices in $G-\{x,y\}$,
there exists a rainbow $x,y$-path in $G$ of length at most $4\ell+4$ that avoids $F$ and  $C$.
\end{lemma}

\begin{proof} 
Since $G$ is a $(d, \lam, \eps)$-expander, $\delta(G) \ge d/2 = d(G)/2$ by Lemma~\ref{lem:d-expand}.
Place each colour used in $G$ in group 1 or 2 with probability $1/2$. For $i \in \{1,2\}$, let $G_i$ denote the spanning subgraph of $G$ that consists of edges whose colour lies in group $i$.
Using Chernoff bounds, we can show that, with positive probability, every $v\in V(G)$ satisfies $d_{G_i}(v)=(\frac{1}{2}\pm o(1)) d_G(v)$.
In particular, when $n$ is sufficiently large, one can ensure that

\medskip

\begin{claim} \label{claim:lower-bound-density}
For $i\in \{1,2\}$, we have $d(G_i)\geq \frac{1}{3} d$ and $\delta(G_i)\geq \frac{1}{3} d(G_i)$.
\end{claim}

By Corollary~\ref{cor:covering}, for $i=1,2$,
$G_i$ contains edge-disjoint subgraphs $G_{i,1}, \dots, G_{i, s_i}$
such that each $j\in [s_i]$, $G_{i,j}$ is a $(d(G_{i,j}),\lambda,\varepsilon)$-expander with $d(G_{i,j})\geq \frac{\varepsilon d(G_i)}{2}\geq  \frac{\varepsilon d}{6}$.
and $e(\cup_{j=1}^{s_i} G_{i,j}) \geq 
(1-\varepsilon) e(G_i)$. 

\medskip

\begin{claim}\label{claim:upperboundonexpanders}
For $i\in\{1,2\}$, we have $s_i\leq \frac{6K}{\varepsilon}$.
\end{claim}

\begin{proof}[Proof of Claim~\ref{claim:upperboundonexpanders}.]
By our discussion, for each $i\in \{1,2\}$ and $j \in \{1,\dots, s_i\}$ we have $d(G_{i,j})\geq \frac{\varepsilon}{6}d$.
By the assumption of the lemma, we must have $v(G_{i,j})>\frac{n}{K}$.  Hence $e(G_{i,j})=\frac{1}{2} v(G_{i,j}) d(G_{i,j})
\geq \frac{1}{2} \frac{n}{K}\frac{\varepsilon d}{6}=\frac{\varepsilon e(G)}{6K}$. This implies that $s_i\leq \frac{6K}{\varepsilon}$. 
\end{proof}
\medskip

Let $S_x$ denote the set of vertices in $G$ which can be reached from $x$ by a rainbow path of length at most $\ell+1$ avoiding $F$ and $C$. Since $G$ is a $(d,\lambda,  \varepsilon)$-expander and $d\geq \frac{4(\ell+2M)}{\lambda}$ and $|F|, |C| \le L \le M$, by Lemma~\ref{lem:reachable-new}, $|S_x|\geq (1-\varepsilon)n$.
Let $J_1=\{j\in [s_1]: V(G_{1,j})\cap S_x\neq \emptyset\}$. 
For each $j\in J_1$, fix a vertex $x_j\in S_x\cap V(G_{1,j})$ and a
rainbow path $P_j$ from $x$ to $x_j$ in $G$ whose length is at most $\ell+1$ avoiding $F$ and $C$. 
Let $F_1 := \cup_{j \in J_1} P_j$. 

 Let $S_y$ to be the set of vertices in $G$ which can be reached from $y$ by a rainbow path of length at most $\ell+1$ avoiding  $F\cup V(F_1)$ and $C \cup C(F_1)$. 
 
 Note that since 
 $$|F\cup V(F_1)|,  |C\cup C(F_1)| \leq  L+ s_1 (\ell+2) \leq  L+ \frac{6K}{\eps} (\ell+2)\leq M  $$
and $d\geq \frac{4(\ell+2M)}{\lambda}$,
by Lemma~\ref{lem:reachable-new},  $|S_y| \ge (1 - \varepsilon)n$.  
Let $J_2=\{j\in [s_2]: V(G_{2,j})\cap S_y\neq \emptyset\}$. 
For each $j\in J_2$, fix a vertex $y_j\in S_y\cap V(G_{2,j})$ and a
rainbow path $Q_j$ from $y$ to $y_j$ in $G$ whose length is at most $\ell+1$
avoiding $F\cup V(F_1)$ and $C \cup C(F_1)$. Let $F_2=\cup_{j\in J_2}Q_j$.
Note that 
\begin{equation} \label{eq:avoid-threshold}
|F\cup V(F_1)\cup V(F_2)|,  |C\cup C(F_1)\cup C(F_2)| \leq  L+(s_1+s_2)(\ell+2)\leq L+ \frac{12 K}{\eps}(\ell+2) \leq M.
\end{equation}

For each $j\in J_1$, let $U_{1,j}$ denote the set of vertices in $G_{1,j}$ that are reachable from $x_j$ by a rainbow path in $G_{1,j}$ of length at most $\left \lceil\frac{4\log v(G_{1,j})}{\lambda}\right \rceil+1 \leq \ell+1$ that avoids $F\cup V(F_1)\cup V(F_2)$ and $C\cup C(F_1)\cup C(F_2)$. 
For each $j\in J_2$, let $U_{2,j}$ denote the set of vertices in $G_{2,j}$ that are reachable from $y_j$ by a rainbow path in $G_{2,j}$ of length at most $\left \lceil\frac{4\log v(G_{2,j})}{\lambda}\right \rceil+1 \leq \ell+1$ that avoids $F\cup V(F_1)\cup V(F_2)$ and $C\cup C(F_1)\cup C(F_2)$. For each $i=1,2$, $j\in J_i$, 
since $G_{i,j}$ is a $(d(G_{i,j}),\lambda, \varepsilon)$-expander and $d(G_{i,j})\geq \frac{\varepsilon d}{6}
\geq\frac{4(\ell+2M)}{\lambda}$, by Lemma~\ref{lem:reachable-new},  we have 
\begin{enumerate}
    \item $|U_{i,j}| \ge (1- \varepsilon) v(G_{i,j})$
    \item $e(G_{i,j}[U_{i,j}])\geq (1-\lambda/2-\varepsilon)e(G_{i,j})\geq (1-2\varepsilon)e(G_{i,j})$.
\end{enumerate}

For $i=1,2$, let $G'_i=\cup_{j\in J_i} G_{i,j}[U_{i,j}]$.  
\begin{claim} \label{claim:G'-size}
For $i=1,2$, we have
$e (G_i')\geq (1-6\varepsilon) e(G_i)$ and $v(G_i') \geq (1-15\varepsilon) n$.
\end{claim}
\begin{proof}
By our discussion above, for each $i\in \{1,2\}, j\in J_i$,  
$e(G_{i,j}[U_{i,j}])\geq (1-2\varepsilon) e(G_{i,j})$.
Since the $G_{i,j}$'s are pairwise edge-disjoint, we have
\begin{equation} \label{eq:Gi-prime}
e(G'_i)=\sum_{j\in J_i} e(G_{i,j}[U_{i,j}])\geq (1-2\varepsilon)e(\cup_{j\in J_i} G_{i,j}).
\end{equation}

By the definition of $J_1$ and $J_2$, $V(\cup_{j \not \in J_1}{G_{1,j}}) \subseteq V(G)\setminus S_x$ and
$V(\cup_{j \not \in J_2}{G_{2,j}}) \subseteq V(G)\setminus S_y$. 
Since $G$ is $d$-minimal with $d=d(G)$, $e(G[V(G)\setminus S_x]) \leq d|V(G)\setminus S_x|/2 \leq d \varepsilon n/2 = \varepsilon e(G)$
and similarly $e(G[V(G)\setminus S_y]) \leq  \varepsilon e(G)$. Hence, for $i=1,2$,
$e(\cup_{j \not \in J_i}{G_{i,j}}) \leq \varepsilon e(G) \leq 3 \varepsilon e(G_i)$ and thus $e(\cup_{j\in J_i}{G_{i,j}})\geq (1-3\varepsilon) e(G_i)$. By \eqref{eq:Gi-prime}, we have for $i=1,2$
 $$ e\left( G_i'\right)\geq (1-2\varepsilon) e(\cup_{j\in J_i}{G_{i,j}}) \geq
 (1-2\varepsilon)(1-3\varepsilon) e(G_i)  \geq  (1-5\varepsilon)e(G_i).$$
 
 By Claim~\ref{claim:lower-bound-density}, $\delta(G_i) \ge \frac{1}{3}d(G_i)$, so by Lemma~\ref{lem:edge-to-vertex} it follows that  $v(G_i') \geq (1-15\varepsilon) v(G_i)$.
 \end{proof}

By Claim~\ref{claim:G'-size} and the fact that $\varepsilon \le \frac{1}{40}$, $V(G_1')\cap V(G_2')\neq \emptyset$.
Let $z$ be any vertex in $V(G'_1)\cap V(G'_2)$. Suppose $z\in U_{1,i}\cap U_{2,j}$.
By definition, there exists a rainbow path $P_i'$ in $G_1$  going from $x_i$ to $z$ and a rainbow path $Q_j'$ in 
$G_2$ going from $y_j$ to $z$, both of which have length at most $\ell+1$ and
avoid $F\cup V(F_1)\cup V(F_2)$ and $C\cup C(F_1)\cup C(F_2)$.  Since $G_1, G_2$ are colour-disjoint,
and $P'_i\subseteq G_1, Q'_j\subseteq G_2$,  $C(P_i')\cap C(Q_j')=\emptyset$.
Also, since $P'_i$ and $Q'_j$ avoid $F\cup F_1\cup F_2$ while $P_i, Q_j\subseteq F_1\cup F_2$, $P'_i$ and $Q'_j$ are colour-disjoint from
$P_i, Q_j$. Finally, recall that $Q_j$ avoids $F_1$ and hence $C(P_i)\cap C(Q_j)=\emptyset$.
Hence, $P_i\cup P_i'\cup Q_j'\cup Q_j$ is a rainbow walk of length at most $4\ell+4$ from $x$ to $y$, that avoids $F$ and $C$ and hence contains a rainbow $x,y$-path of length at most $4\ell+4$ that avoids $F$ and $C$, as desired.
\end{proof}

\begin{corollary} \label{cor:rainbow-subdivision}
For every integer $t \ge 2$, $0<\varepsilon \le \frac{1}{40}$ and sufficiently large $n$ the following holds.  
Let $\lambda = \frac{\varepsilon}{2\log n}$. Let $\ell=\left \lceil\frac{4\log n}{\lambda}\right \rceil$. Let $L=\binom{t}{2}(4\ell+4)+t$.
Let $K$ be a positive integer. Let $M \ge L+ \frac{12 K}{\eps}(\ell+2)$.
Let $d$ be a real that satisfies  $d\geq \frac{48(\ell+2M)}{\varepsilon\lambda}$.
Let $G$ be a properly edge-coloured $n$-vertex graph with $d(G)\geq d$.
Suppose $G$ contains no subgraph $G'$ satisfying both $d(G')\geq \frac{\varepsilon d(G)}{12}$ and $v(G')\leq \frac{n}{K}$. 
Then $G$ contains a rainbow subdivision of $K_t$, where each edge of $K_t$ is replaced with a path of length at most $4\ell+4$.
\end{corollary}

\begin{proof}  
By Lemma~\ref{lem:existence-expanders-new}, $G$ contains a subgraph $H$ which is a $(d(H),\lambda, \varepsilon)$-expander with $d(H)\geq \frac{d(G)}{2}$. 
We may assume that $H$ does not have a subgraph $H'$ on at most $\frac{v(H)}{K}$ vertices with $d(H') \ge \frac{\eps d(H)}{6}$. Indeed, otherwise, such a subgraph $H'$ would be a subgraph
of $G$ on at most $\frac{n}{K}$ vertices with $d(H')\geq \frac{\varepsilon d(H)}{6} \geq \frac{\varepsilon d(G)}{12}$, contradicting our assumption about $G$.

Let $S=\{x_1,\dots, x_t\}$ be any set of $t$ distinct vertices in $H$.
Let $A$ be a maximal collection of pairs $(i, j)$ with $1 \le i < j \le t$ such that there exist paths $P_{i, j}$ for $(i,j)\in A$ satisfying the following. 
\begin{enumerate}
    \item For each $(i,j)\in A$,
$P_{i, j}$ is  an $x_i, x_j$-path of length at most $4\ell+4$ in $H$ such that $V(P_{i,j})\cap \{x_1,\dots, x_t\}=\{x_i, x_j\}$.
\item The sets $V(P_{i,j})\setminus\{x_i,x_j\}$ are pairwise disjoint over different $(i,j)\in A$
and $\cup_{(i,j)\in A} P_{i,j}$ is rainbow. 
\end{enumerate}
If $A$ consists of all pairs $(i, j)$ with $1 \le i < j \le t$, then the union of the paths $P_{i, j}$ forms the required subdivision of $K_t$ and we are done. Hence, we may assume that there exist $i,j$
such that $1 \le i < j \le t$ and $(i, j) \notin A$. 

Let $C$ be the set of colours used on  $\cup_{(k,\ell)\in A} P_{k, \ell}$ and
let $F=V(\cup_{(k,\ell)\in A} P_{k,\ell})\cup S\setminus \{x_i,x_j\}$. Then $|C|,|F|\leq \binom{t}{2}(4\ell+4)+t=L$. Since $n$ is sufficiently large and $d(H)\geq \frac{d}{2}$, we may assume $v(H)$ is sufficiently large. Also, $d(H)\geq \frac{d}{2}\geq \frac{24(\ell+2M)}{\lambda}$.
Hence, by Lemma \ref{rainbow-connection}, $H$ contains
a rainbow $x_i, x_j$-path $P_{i, j}$ that avoids $C$ and $F$ and whose length is at most $4\ell+4$. The pair $(i, j)$ could thus be added to $A$, contradicting maximality of $A$. 
\end{proof}

Now we are ready to prove Theorem~\ref{thm:main}.

\medskip

\textbf{Proof of Theorem \ref{thm:main}:} 
Let $\varepsilon=\frac{1}{40}$. The choice of $n_0$ is not explicit; we assume it is sufficiently large with respect to $t$ so that various inequalities below hold.

Choose $c= \max\{\sqrt{\log{n_0}}, 2\log(\frac{12}{\varepsilon})\}$. By our choice of $c$, $e^{\frac{c}{2}}\geq \frac{12}{\varepsilon}$.
We may assume $n\geq n_0$. 
Let $\lambda = \frac{\varepsilon}{2\log n}$. 
Let $\ell=\left \lceil\frac{4\log n}{\lambda}\right \rceil$. Let $L=\binom{t}{2}(4\ell+4)+t$.
Let $K=\left \lceil e^{\sqrt{\log n}}\right \rceil$. 
Let $M=L+ \frac{12 K}{\eps}(\ell+2)$.
Note that by our choice of $c$ and $n_0$, we have $e^{\frac{c}{2}\sqrt{\log n}}\geq \frac{48(\ell+2M)}{\varepsilon\lambda}$.

Let $G$ be a properly edge-coloured graph on $n$ vertices with at least $n e^{c\sqrt{\log n}}$ edges, so $d(G) \ge e^{c\sqrt{\log n}}$. 
Let $G_0=G$. We define a sequence of subgraphs of $G$ as follows. For each $i\geq 0$,
suppose $G_0,\dots, G_i$ have been defined. If $G_i$ contains a subgraph $F$ with
$d(F)\geq \frac{\varepsilon}{12} d(G_i)$ and $v(F)\leq \frac{v(G_i)}{K}$, let $G_{i+1}=F$.
If no such subgraph of $G_i$ exists then we stop the process. Suppose the process stops when considering the graph $G_m$, so no such subgraph $F$ can be found inside $G_m$. Then $1\leq v(G_m)\leq \frac{n}{K^m}$. Hence $m\leq \log_K n \leq \frac{\log n}{\log e^{\sqrt{\log n}}}=\sqrt{\log n}$. Also, since $e^\frac{c}{2}\geq \frac{12}{\varepsilon}$ by our assumption, we have
\[d(G_m)\geq \left(\frac{\varepsilon}{12}\right)^m d(G)
\geq \left(\frac{\varepsilon}{12}\right)^{\sqrt{\log n}} e^{c\sqrt{\log n}} \geq
e^{\frac{c}{2}\sqrt{\log n}}.\]

By our choice of $c$ and $n_0$, we have $d(G_m)\geq e^{\frac{c}{2}\sqrt{\log n}}
\geq \frac{48(\ell+2M)}{\varepsilon\lambda} \geq  \frac{48(\ell'+2M')}{\varepsilon\lambda'}$, where $\lambda'= \frac{\eps}{2 \log v(G_m)}, \ell'= \left \lceil \frac{8 \log v(G_m)}{\lam} \right \rceil $, $L'=\binom{t}{2}(4\ell'+4)+t$, $M'=L'+ \frac{12 K}{\eps}(\ell'+2).$
Also, $v(G_m)>d(G_m)\geq e^{\frac{c}{2}\sqrt{\log n}}>n_0$,
 and by the definition of $G_m$, $G_m$ contains no subgraph $F$ satisfying
both $d(F)\geq \frac{\eps}{12} d(G_m)$ and $v(F)\leq \frac{v(G_m)}{K}$.

So by Corollary \ref{cor:rainbow-subdivision}, $G_m$ contains 
a rainbow subdivision of $K_t$, where each edge of $K_t$ is replaced with
a path of length at most $4\ell+4\leq 1300\log^2 n$.
This completes our proof.

\section{Concluding remarks}

This paper establishes that a properly edge-coloured graph on $n$ vertices with average degree at least $(e^{c\sqrt{\log{n}}})$ contains a  rainbow subdivision of a clique. The proof of Theorem~\ref{thm:main} in fact
explicitly gives the following slightly stronger result.  Given a positive integer $t\geq 2$, we call a properly edge-coloured graph $G$ {\it rainbow $t$-linkable} if for any set of $t$ vertices $x_1,\dots, x_t$ in $G$ there is a rainbow
$K_t$-subdivision with $x_1,\dots, x_t$ playing the role of the $t$ branching vertices. For every $t\geq 2$, there is a constant $c$ such that for every $n \ge t$, if $G$ is a properly coloured $n$-vertex graph with average degree at least $e^{c\sqrt{\log n}}$, then $G$ contains a subgraph $H$ with average degree at least $\Omega(e^{\frac{c}{2}\sqrt{\log n}})$ that is rainbow $t$-linkable.  

After the initial submission of our paper, there has been further improvements on the average degree that enforces the existence of a rainbow cycle and a rainbow subdivision of a clique. For rainbow cycles, Janzer's result was improved first by Tomon~\cite{tomon} to  $(\log n)^{2+o(1)}$ and by Kim, Lee, Liu, Tran~\cite{joonkyung} to $100(\log{n})^2$.  As for the family of subdivisions of a clique, in a subsequent paper with Letzter~\cite{4-person-paper} we improved the average degree bound to $(\log{n})^{60}$. Subsequently, Tomon~\cite{tomon} improved this average degree lower bound to  $(\log{n})^{6+o(1)}$, and finally Wang~\cite{Wang} improved it to  $(\log{n})^{2+o(1)}$. We pose the following question (also stated in~\cite{4-person-paper}).

\begin{question} 
    Given $t\geq 3$, what is the smallest $c$ such that for all sufficiently large $n$, if $G$ is a properly edge-coloured graph on $n$ vertices with $\Omega(n(\log{n})^c)$  many edges then $G$ contains a rainbow subdivision of $K_t$?  In particular, is $c=1$?
\end{question}

In the above question, our guess of $c=1$ is motivated by the only lower bound of order $\Omega(n\log{n})$ given by the construction of a colouring of the hypercube by Keevash, Mubayi, Sudakov and Verstra\"ete \cite{keevash2007rainbow} described earlier in the introduction.
 
\section{Acknowledgments}

We thank Shoham Letzter for carefully reading an early draft of this paper and for helpful comments. We also thank the referees for their valuable feedback which improved the presentation of the paper significantly.

\bibliographystyle{amsabbrv}
\bibliography{ref}

\end{document}